\documentclass{article}

\usepackage[final]{neurips_2019}

\setcitestyle{numbers,square,comma}

\usepackage[utf8]{inputenc} 
\usepackage[T1]{fontenc}    
\usepackage{url}            

\usepackage{microtype}
\usepackage{graphicx}
\usepackage{subfigure}
\usepackage{booktabs} 
\usepackage{hyperref}
\usepackage{longtable}
\usepackage{amsfonts,amsmath,amssymb,graphicx,color}

\usepackage{psfrag}
\usepackage{epstopdf}
\usepackage{algorithm, algorithmic}

\usepackage{amsthm}
\usepackage{psfrag}
\usepackage{pdflscape}
\usepackage{multirow}
\usepackage{blindtext} 
\usepackage{enumitem}
\hypersetup{colorlinks=true}
\usepackage{framed}
\usepackage{footnote}
\usepackage[dvipsnames]{xcolor}
\usepackage{soul}
\usepackage{multirow}
\usepackage{grffile}
\usepackage{nicefrac}

\usepackage{chngcntr}

\newtheorem{thm}{Theorem}[section]

\newtheorem{assum}[thm]{Assumptions}

\counterwithout{equation}{section}

\DeclareMathOperator{\Var}{Var}

\newcommand{\defeq}{:=}

\newcommand{\R}{{\mathbb R}}        

\newcommand{\longonly}[1]{ } 

\newcommand{\FD}{\textrm{FD}} 

\def\noprint#1{}

\def\E{\mathbb{E}}

\newcommand{\CASE}[1]{\STATE \textbf{case} #1\textbf{:} \begin{ALC@g}}
	\newcommand{\ENDCASE}{\end{ALC@g}}

\newcommand{\DEFAULT}{\STATE \textbf{default:} \begin{ALC@g}}
	\newcommand{\ENDDEFAULT}{\end{ALC@g}}
\newcommand{\DEFAULTLINE}[1]{\STATE \textbf{default:} }

\title{Adaptive Sampling Quasi-Newton Methods for Derivative-Free Stochastic Optimization}

\author{%
  Raghu Bollapragada\\
  Argonne National Laboratory\\
  Lemont, IL, USA\\
  \texttt{raghu.bollapragada@utexas.edu} \\
  \And
  Stefan M.\ Wild\\
  Argonne National Laboratory\\
  Lemont, IL, USA\\
  \texttt{wild@anl.gov} 
}

\begin{document}

\maketitle

\begin{abstract}
We consider stochastic zero-order optimization problems, which 
arise in settings from simulation optimization to reinforcement learning.
We propose an adaptive sampling quasi-Newton method where we 
estimate the gradients of a stochastic function using finite differences
within a common random number framework. 
We employ modified versions of 
a 
\emph{norm test} and 
an
\emph{inner product quasi-Newton test} to control the sample sizes used in the stochastic approximations.  We provide preliminary numerical experiments to illustrate potential performance benefits of the proposed method.
\end{abstract}

\section{Introduction}
We consider unconstrained stochastic optimization problems of the form 
\begin{equation}
\label{eq:prob}
\min_{x \in \Re^d} F(x) = E_{\zeta}\left[f(x,\zeta)\right],
\end{equation}
where one has access only to a zero-order oracle (i.e., a black-box procedure that outputs realizations of the stochastic function values $f(x,\zeta)$ and cannot access explicit estimates of the gradient 
$\nabla_x f(x,\zeta)$). 
Such stochastic optimization problems arise in a plethora of science and engineering applications, from simulation optimization \cite{BCMS2018,Fu2005,Kim2014,Pasupathy2013} 
to reinforcement learning \cite{Bertsimas2019b,Mania2018,Salimans2017}.
Several
methods have been proposed 
to solve 
such ``derivative-free'' 
problems \cite{AudetHare2017,LMW2019AN}.

We propose 
finite-difference stochastic quasi-Newton methods for solving \eqref{eq:prob} 
by exploiting common random number (CRN) evaluations of $f$.
The CRN setting allows us to define subsampled gradient estimators
\begin{eqnarray}
\label{eq:FFD}
\left[\nabla F^{\FD}_i(x)\right]_j 
&\defeq& \frac{f(x + \nu e_j, \zeta_i) - f(x, \zeta_i)}{\nu}, 
\; j=1,\ldots,d
\\ 
\label{eq: batch_FD}
\nabla F^{\FD}_{S_k}(x) &\defeq& \frac{1}{|S_k|}\sum_{\zeta_i \in S_k}\nabla F_i^{\FD}(x),
\end{eqnarray}
which employ forward differences for the i.i.d.\ samples of $\zeta$ in the set $S_k$ and whereby the parameter $\nu$ needs to account only for numerical errors.
%
This gradient estimation has two sources of error: error due to the finite-difference approximation and error due to the stochastic approximation. The error due to stochastic approximation depends on the number of samples $|S_k|$ used in the estimation. Using too few samples affects the stability of a method using the estimates; using a large number of samples results in computational inefficiency. 
For settings where gradient information is available, researchers have developed practical tests to adaptively increase the sample sizes used in the stochastic approximations and have supported these tests with  global convergence results  \cite{bollapragada2018adaptive, pmlr-v80-bollapragada18a, byrd2012sample,pasupathy2018sampling}.
In this paper we modify these tests to address the challenges associated with the finite-difference approximation errors, 
and we demonstrate the resulting method on simple test problems.


\section{A Derivative-Free Stochastic Quasi-Newton Algorithm}
\label{sec:alg}
The update form of a finite-difference, derivative-free stochastic quasi-Newton method is given by  
\begin{equation}
\label{eq:iter}
x_{k+1} = x_k - \alpha_k H_k \nabla F_{S_k}^{\rm FD}(x_k),
\end{equation}
where $\alpha_k > 0$ is the steplength, $H_k$ is a positive-definite quasi-Newton matrix, and $ \nabla F^{\FD}_{S_k}(x_k)$ is the (e.g., forward) finite-difference subsampled (or batch) gradient estimate defined by \eqref{eq: batch_FD}. 
We propose to control the sample sizes $|S_k|$ over the course of optimization 
to achieve fast convergence by using two different strategies adapted to the setting where no gradient information is available.

\paragraph{Norm test.} The first test is based on the norm condition \cite{byrd2012sample,hashemi2014adaptive}:
\begin{equation}\label{eq:ideal_norm}
\E_{S_k}[\|\nabla F^{\FD}_{S_k}(x_k) - \nabla F(x_k)\|^2] \leq \theta^2 \|\nabla F(x_k)\|^2,
\end{equation}
which is used for controlling the sample sizes in subsampled gradient methods. We note that the finite-difference approximation error (first pair of terms) in $\nabla F^{\FD}_{S_k}(x_k) - \nabla F(x_k)=\nabla F^{\FD}_{S_k}(x_k) - \nabla F_{S_k}(x_k)+ \nabla F_{S_k}(x_k) - \nabla F(x_k)$ is nonzero and can be upper bounded (independent of $S_k$) for any function $f(\cdot,\zeta)$ with $L$-Lipschitz continuous gradients:
\begin{align}
\|\nabla F^{\FD}_{S_k}(x_k) - \nabla F_{S_k}(x_k)\|^2  &\leq \left(\tfrac{L\nu \sqrt{d}}{2}\right)^2,
\qquad \mbox{where }  \nabla F_{S_k}(x_k) := \tfrac{\sum_{\zeta_i \in S_k} \nabla_x f(x_k,\zeta_i)}{|S_k|};\label{eq:fdbound}
\end{align}
the proof is given in the supplementary material. Therefore, one cannot always satisfy \eqref{eq:ideal_norm}; moreover, satisfying \eqref{eq:ideal_norm} might be too restrictive. 
Instead, we propose to look at the norm condition based on the finite-difference subsampled gradient estimators. That is, we use the condition
\begin{equation}
\label{eq:ideal_normFD}
\E_{S_k}[\|\nabla F^{\FD}_{S_k}(x_k) - \nabla F^{\FD}(x_k)\|^2] \leq \theta^2 \|\nabla F^{\FD}(x_k)\|^2,  \quad \nabla F^{\FD}(x_k) \defeq \left[\tfrac{F(x_k + \nu e_j) - F(x_k)}{\nu}\right]_{j=1}^{d},
\end{equation} 
for which $\E_{S_k}[\nabla F^{\FD}_{S_k}(x_k)] = \nabla F^{\FD}(x_k)$ and which corresponds to a norm condition where the right-hand side of \eqref{eq:ideal_norm} is relaxed. \longonly{(see \eqref{eq:relax_norm} in the supplementary material). }
The left-hand side of \eqref{eq:ideal_normFD} is difficult to compute but can be bounded by the true variance of individual finite-difference gradient estimators; this results in 
\begin{equation}
\label{eq:popl_normFD}
\frac{\E_{i}[\|\nabla F^{\FD}_{i}(x_k) - \nabla F^{\FD}(x_k)\|^2]}{|S_k|} \leq \theta^2 \|\nabla F^{\FD}(x_k)\|^2.
\end{equation}  
Approximating the true expected gradient and variance with sample gradient and variance estimates, respectively, yields the  \emph{practical finite-difference norm test}:
\begin{equation}\label{eq:sample_normFD}
\frac{\tfrac{1}{|S_k^v| - 1}\sum_{\zeta_i \in S_k^v} \|\nabla F_i^{\FD}(x_k) - \nabla F_{S_k}^{\FD}(x_k)\|^2}{|S_k|} \leq \theta^2  \|\nabla F_{S_k}^{\FD}(x_k)\|^2,
\end{equation}
where $S_k^v \subseteq S_k$ is a subset of the current sample (batch). 
In our algorithm, we test condition \eqref{eq:sample_normFD};  whenever it is not satisfied, we control $|S_k|$ such that the condition is satisfied.

\paragraph{Inner product quasi-Newton test.} 
The norm condition controls the variance in the gradient estimation and does not utilize observed quasi-Newton information to control the sample sizes. Recently, Bollapragada et al.\ \cite{pmlr-v80-bollapragada18a} 
proposed to control the sample sizes used in the gradient estimation by ensuring that the stochastic quasi-Newton directions make an acute angle with the true quasi-Newton direction with high probability. That is,
\begin{equation} \label{eq:ip}
(H_k \nabla F_{S_k}^{\FD}(x_k))^T(H_k \nabla F(x_k)) > 0
\end{equation}
holds with high probability. This condition can be satisfied in expectation at points that are sufficiently far away from the stationary points; that is, for points $x_k$ such that $\|\nabla F(x_k)\| > \tfrac{\lambda_{\max}(H_k) L\nu\sqrt{d}}{2\lambda_{\min}(H_k)}$, where $\lambda_{\max}(H_k)$ and $\lambda_{\min}(H_k) > 0$, are the largest and smallest eigenvalues of $H_k$, respectively (see supplementary material). Hence, the condition \eqref{eq:ip} can be satisfied with high probability at points $x_k$ that are sufficiently far away from being stationary points. We must control the variance in the left-hand side of \eqref{eq:ip} 
to achieve this objective. We note that the quantity $H_k\nabla F(x_k)$ cannot be computed directly; however, it can be approximated by $H_k\nabla F^{\FD}(x_k)$. The condition is given as 
\begin{align}\label{eq:ideal_ipqnFD}
\E_{S_k}\left[\left( (H_k \nabla F_{S_k}^{\FD}(x_k))^T(H_k \nabla F^{\FD}(x_k)) - \left\|H_k \nabla F^{\FD}(x_k)\right\|^2 \right)^2\right]
&\leq \theta^2 \|H_k \nabla F^{\FD}(x_k)\|^4,
\end{align}
where $\E_{S_k}\left[H_k \nabla F_{S_k}^{\FD}(x_k)\right]=H_k\nabla F^{\FD}(x_k)$. The left-hand side of \eqref{eq:ideal_ipqnFD} can be bounded by the true variance, as done before.  Therefore, the following condition is sufficient for ensuring that \eqref{eq:ideal_ipqnFD} holds:
\begin{align}
\frac{\E_{i}\left[\left((H_k \nabla F_{i}^{\FD}(x_k))^T(H_k \nabla F^{\FD}(x_k)) - \|H_k \nabla F^{\FD}(x_k)\|^2\right)^2\right]}{|S_k|}
&\leq \theta^2 \|H_k \nabla F^{\FD}(x_k)\|^4. \label{eq:popl_ipqnFD}
\end{align}
Approximating the true expected gradient and variance with sample gradient and variance estimates results in the \emph{practical finite-difference inner product quasi-Newton test:}
\begin{equation}\label{eq:sample_ipqnFD}
\frac{\tfrac{1}{|S_k^v| - 1} \sum_{\zeta_i \in S_k^v} \left((H_k\nabla F_{S_k}^{\FD}(x_k))^TH_k\nabla F_{i}^{\FD}(x_k) - \|H_k\nabla F^{\FD}_{S_k}\|^2\right)^2}{|S_k|} \leq \theta^2  \|H_k\nabla F_{S_k}^{\FD}(x_k)\|^4 ,
\end{equation}
where $S_k^v \subseteq S_k$ is a subset of the current sample (batch).
This variance computation requires only one additional Hessian-vector product (i.e., the product of $H_k$ with $H_k \nabla F_{S_k}^{\FD}(x_k)$).
In our algorithm, we test the condition \eqref{eq:sample_ipqnFD}; whenever it is not satisfied, we control $|S_k|$ to satisfy the condition. 

\paragraph{Finite-difference parameter and steplength selection.}
We select the finite-difference parameter $\nu$ by minimizing an upper bound on the error in the gradient approximation. 
Assuming that numerical errors in computing $f(x, \zeta)$ are uniformly bounded by $\epsilon_m$ yields the parameter value
$
\nu^* \defeq 2\sqrt{\tfrac{\epsilon_m}{L}}.
$
\longonly{(see supplementary material).}
This 
finite-difference parameter is 
analogous to
the one 
derived in \cite{more2011edn},
which depends on the variance in stochastic noise;
however, in the CRN setting we need to account only for numerical errors.

We employ a stochastic line search to choose the steplength $\alpha_k$ by using a sufficient decrease condition based on the sampled function values. 
In particular, 
we would like 
$\alpha_k$ to satisfy
\begin{equation}\label{eq:suff_decrease}
F_{S_k}\left(x_k - \alpha_k H_k \nabla F_{S_k}^{\FD} (x_k)\right) \leq F_{S_k}(x_k) - c_1 \alpha_k (\nabla F_{S_k}^{\FD}(x_k))^TH_k\nabla F_{S_k}^{\FD}(x_k),
\end{equation}
where $F_{S_k}(x_k) = \tfrac{1}{|S_k|}\sum_{\zeta_i \in S_k} f(x_k,\zeta_i)$ and $c_1 \in (0,1)$ is a user-specified parameter. We employ a backtracking procedure wherein a trial $\alpha_k$ that does not satisfy \eqref{eq:suff_decrease} is reduced by a fixed fraction $\tau < 1$ (i.e., $\alpha_k \leftarrow \tau \alpha_k$).  
One cannot always satisfy \eqref{eq:suff_decrease}, since the quasi-Newton direction $-H_k \nabla F_{S_k}^{\FD}(x_k)$ may not be a descent direction for the sampled function $F_{S_k}$ at $x_k$. Intuitively, at points where the error in the sample gradient estimation error $\nabla F_{S_k}^{\FD}(x_k) - \nabla F_{S_k}(x_k)$ dominates the sample gradient itself, measured in terms of the matrix norm $\|\cdot\|_{H_k}$, this condition may not be satisfied, and the line search fails. 
In practice, stochastic line search failure is an indication that the method has converged to a neighborhood of a solution for the sampled function $F_{S_k}$, and the solution cannot be further improved. Therefore, one can use the line search failure as an early stopping criterion.

We also note that because of the stochasticity in the function values, it is not guaranteed that a decrease in stochastic function realizations can ensure decrease in the true function. A conservative strategy to address this issue is to choose the initial trial steplength to be small enough such that the increase in function values (when the stochastic approximations are not good) is controlled. Following the strategy proposed in \cite{pmlr-v80-bollapragada18a}, we derive a heuristic to choose the initial steplength as 
\begin{equation}
\label{eq:stepinitial}
\hat \alpha_k = \left(1 + \frac{\Var_{\zeta_i \in S_k^v}\left(\nabla F_i^{FD}(x_k)\right)}{|S_k|\|\nabla F_{S_k}^{\FD}(x_k)\|^2}\right)^{-1},
\end{equation}
where $\Var_{\zeta_i \in S_k^v}\left(\nabla F_i^{FD}(x_k)\right)\defeq \tfrac{\sum_{i \in S_k^v} \|\nabla F_i^{\FD}(x_k) - \nabla F_{S_k}^{\FD}(x_k)\|^2}{|S_k^v| - 1}$ is the sample variance used in \eqref{eq:sample_normFD}. 
\longonly{The reasoning for this choice is provided in the supplementary material.}

\paragraph{Stable quasi-Newton update.}
In BFGS and L-BFGS methods, the inverse Hessian approximation is updated by using the formulae
\begin{equation}
\begin{aligned}
H_{k+1} & = V_k^TH_kV_k + \rho_ks_ks_k^T, \qquad 
\rho_k  = (y_k^Ts_k)^{-1}, \qquad 
V_k  = I - \rho_ky_ks_k^T,
\end{aligned}
\end{equation}
where $s_k = x_{k+1} - x_k$ and $y_k$ is  the difference in the gradients at $x_{k+1}$ and $x_k$. In stochastic settings, several recent works \cite{berahas2016multi,pmlr-v80-bollapragada18a, schraudolph2007stochastic} define $y_k$ as the difference in gradients measured on the same sample $S_k$ to ensure stability in the quasi-Newton approximation. We follow the same approach and define
\begin{equation}\label{full-overlap}
y_k \defeq \nabla F_{S_k}^{\FD}(x_{k+1}) - \nabla F_{S_k}^{\FD}(x_k).
\end{equation}
However, even though computing gradient differences on common sample sets can improve stability, the curvature pair $(y_k,s_k)$ still may not satisfy the condition $y_k^Ts_k > 0$
required to ensure positive definiteness of the quasi-Newton matrix. \longonly{(see supplementary material).} \longonly{In particular, for any $\mu$-strongly convex function $F_{S_k}$, we have that
\begin{align}
y_k^Ts_k  &= (\nabla F_{S_k}^{\FD}(x_{k+1}) - \nabla F_{S_k}^{\FD}(x_k))^Ts_k \nonumber\\
&= (\nabla F_{S_k}(x_{k+1}) - \nabla F_{S_k}(x_k))^Ts_k \nonumber \\
&\quad + (\nabla F_{S_k}^{\FD}(x_{k+1}) - \nabla F_{S_k}(x_{k+1}) + \nabla F_{S_k}(x_k) - \nabla F_{S_k}^{\FD}(x_k))^Ts_k \nonumber \\
&\geq \mu\|s_k\|^2 - (\|\nabla F_{S_k}^{\FD}(x_{k+1}) - \nabla F_{S_k}(x_{k+1})\| + \|\nabla F_{S_k}(x_k) - \nabla F_{S_k}^{\FD}(x_k)\|)\|s_k\| \nonumber \\
&\geq \mu\|s_k\|^2 - L\nu\sqrt{d}\|s_k\|  =  \|s_k\| (\mu \|s_k\| -  L\nu\sqrt{d}) \nonumber
\end{align}
Therefore, the condition $y_k^Ts_k > 0$ is guaranteed to be satisfied when $\|s_k\| > \frac{L\nu\sqrt{d}}{\mu}$. Recently, Xie et al. \cite{xie2019analysis} have proposed to modify the curvature pair update whenever the step $s_k$ is too small so that $y_k^Ts_k > 0$. However, this modification requires knowledge of some unknown problem parameters and may not provide guarantees in the nonconvex case.} \longonly{Therefore, we skip the quasi-Newton update if the following curvature condition is not satisfied:
\begin{equation}
y_k^T s_k > \beta \|s_k\|^2, \quad\mbox{with} \quad  \beta=10^{-2}. 
\end{equation}}
Therefore, in our tests we skip the quasi-Newton update if the condition $y_k^T s_k > \beta \|s_k\|^2$ is not satisfied with $\beta=10^{-2}$.

\paragraph{Convergence results.} We leave to future work the settings under which one can establish global convergence results to a neighborhood of an optimal solution for the proposed methods.

\section{Numerical Experiments}
As a demonstration, we conducted preliminary numerical experiments on stochastic nonlinear least-squares problems
based on a mapping $\phi:\R^d \rightarrow \R^p$ affected by
two forms of stochastic noise:
\begin{align}
\label{eq:nonlinearls}
f_{\rm rel}(x,\zeta)\defeq \frac{1}{1 + \sigma^2} \sum_{j=1}^{p}\phi^2_j(x) (1 + \zeta_j)^2 
\qquad \mbox{and} \qquad 
f_{\rm abs}(x,\zeta) \defeq \sum_{j=1}^{p}(\phi_j(x) + \zeta_j)^2 - \sigma^2,
\end{align}
where $\zeta \sim \mathcal{N}(0, \sigma^2I_{p})$ and $\E_{\zeta}[f(x,\zeta)] = \sum_{j=1}^{p}\phi^2_j(x)$. 
In both cases, the function $f(\cdot,\zeta)$ and the expected function $E_{\zeta}[f(\cdot,\zeta)]$ are twice continuously differentiable.
Here we report results only for $\phi$ defined by the Chebyquad function from \cite{Fletcherdfo65} with $d=30$, $p=45$, and an approximate noise-free value $F^*=0.0174$.
\longonly{We computed the minimum function values $F^*$ by running the L-BFGS method on the noise-free (i.e., $\sigma=0$) problems until $\|\nabla F(x)\|_{\infty} \leq 10^{-10}$ 
	We implemented all the algorithms in MATLAB R2019a and run the experiments on a 64-bit mac machine (machine precision $\epsilon_m = 10^{-16}$) with Intel Core i5@2.4GHz and 8GB RAM.
}

We implemented two variants, ``FD-Norm'' and ``FD-IPQN'', of the proposed algorithm using L-BFGS
with\longonly{ the sample size} $|S_k|$ chosen based on either\longonly{ the \emph{finite-difference norm test} in} \eqref{eq:sample_normFD} or\longonly{ the \emph{inner product quasi-Newton test} given in}
\eqref{eq:sample_ipqnFD}, respectively. 
%
%
We also implemented a \emph{finite-difference stochastic gradient} method (``FD-SG'')\longonly{ where the update is given by} 
$x_{k+1} = x_k - \alpha \nabla_{S_k}^{FD}(x_k)$, 
where $\nabla F_{S_k}^{FD}(x_k)$ is defined in \eqref{eq: batch_FD} and $|S_k| = |S_0|$. We report results for the best version of FD-SG based on tuning the \longonly{constant sample size $|S_0|$ and} constant steplength for each problem (e.g., considering $\alpha_0 = 2^j$, for $j \in \{-20, -19, \ldots, 9, 10\}$).
%
We use the initial sample $|S_0| = 64$ for larger variance ($\sigma = 10^{-3}$) problems and $|S_0| = 2$ for smaller variance problems ($\sigma = 10^{-5}$) in all the methods. 
\longonly{The initial starting point is chosen as $x_0 = 10 x_s$ where $x_s$ is the standard starting point for these problems given in \cite{?}. In the experiments, we report, \emph{total function evaluations}, which represent total number of individual functions evaluated in the gradient computations, curvature pair updates, and stochastic backtracking line search. } 

\begin{figure}[!tb]
	\begin{centering}
		\includegraphics[width=0.40\linewidth]{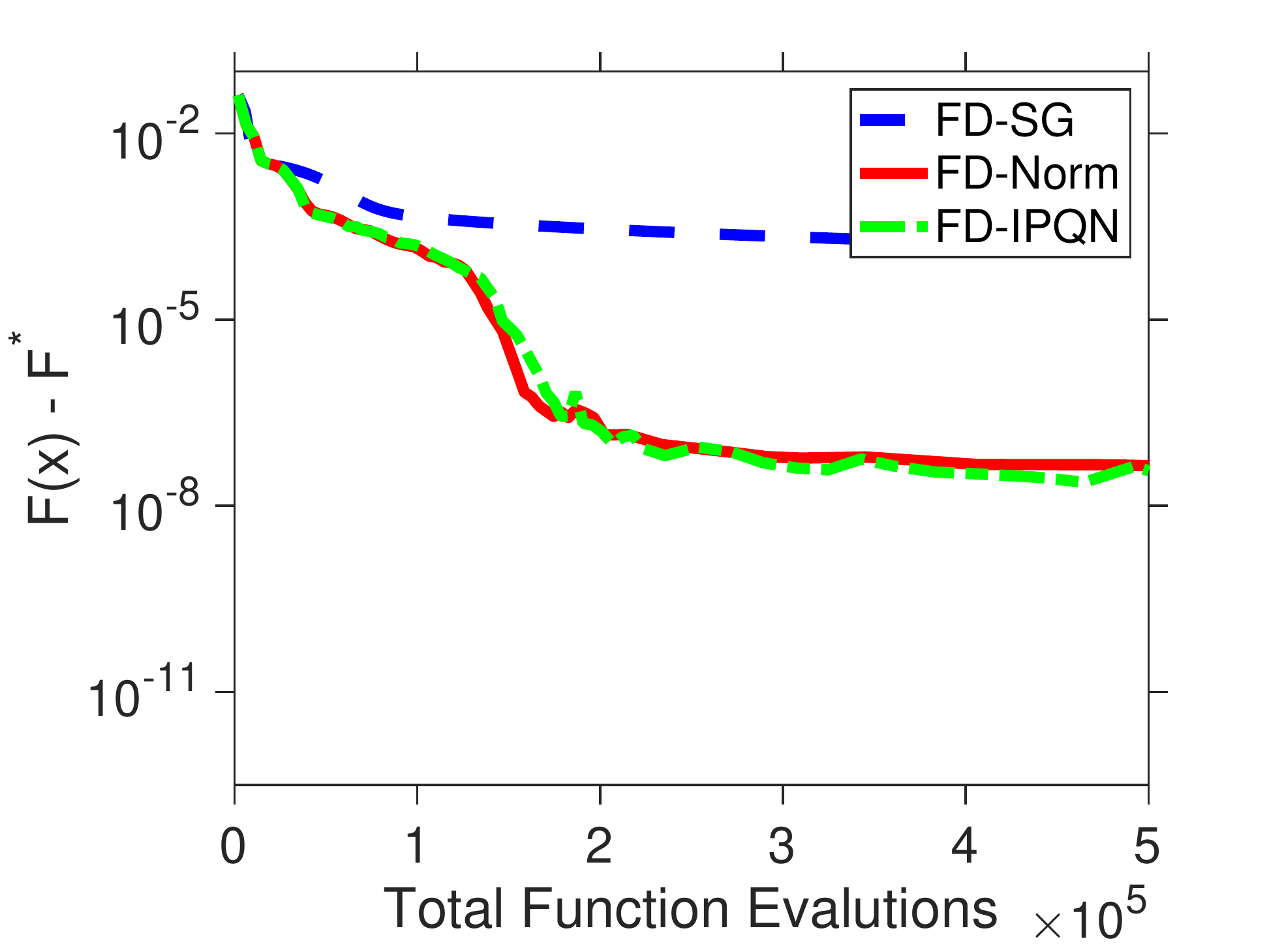} 
		\includegraphics[width=0.40\linewidth]{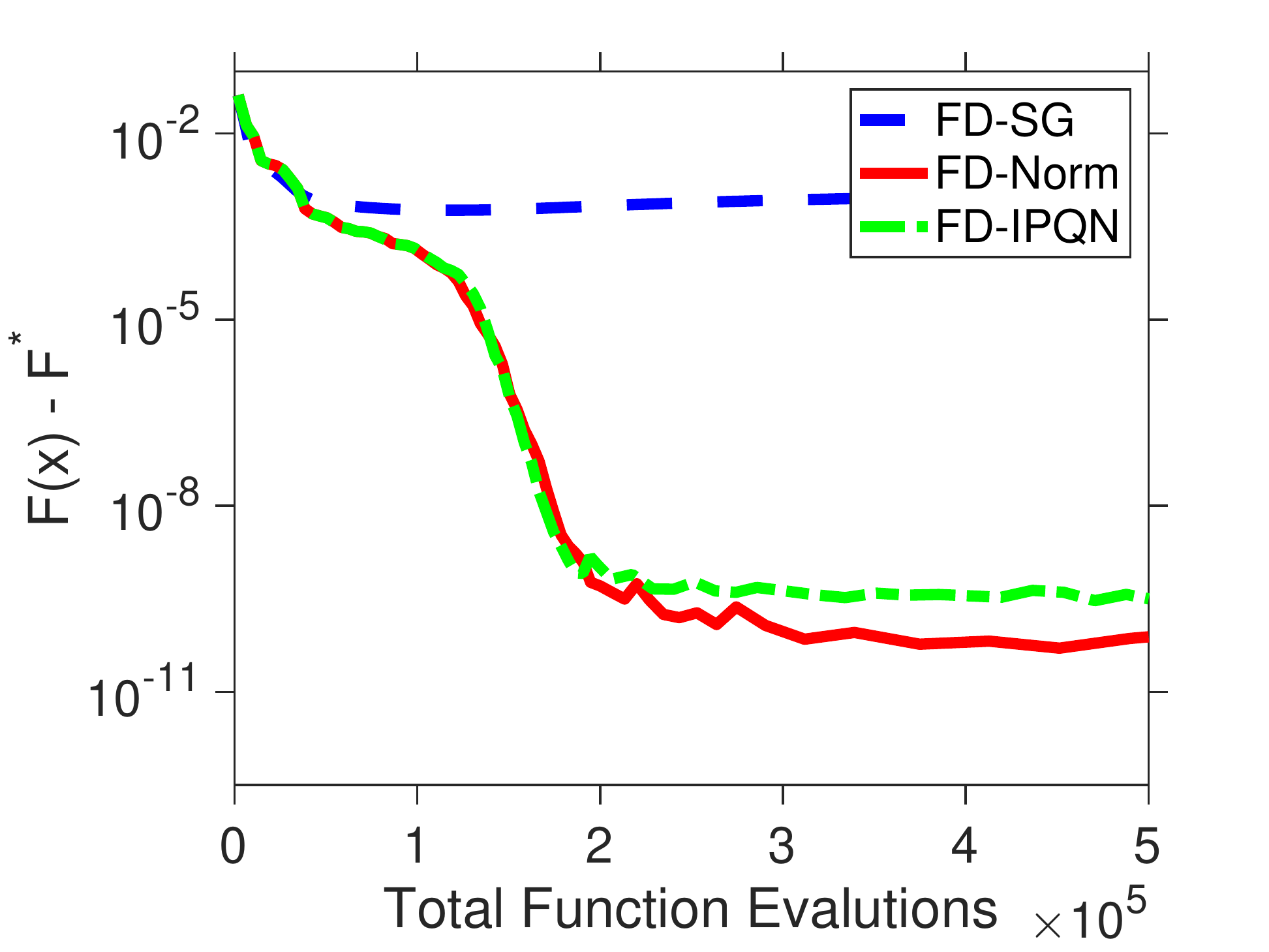}\\
		\includegraphics[width=0.40\linewidth]{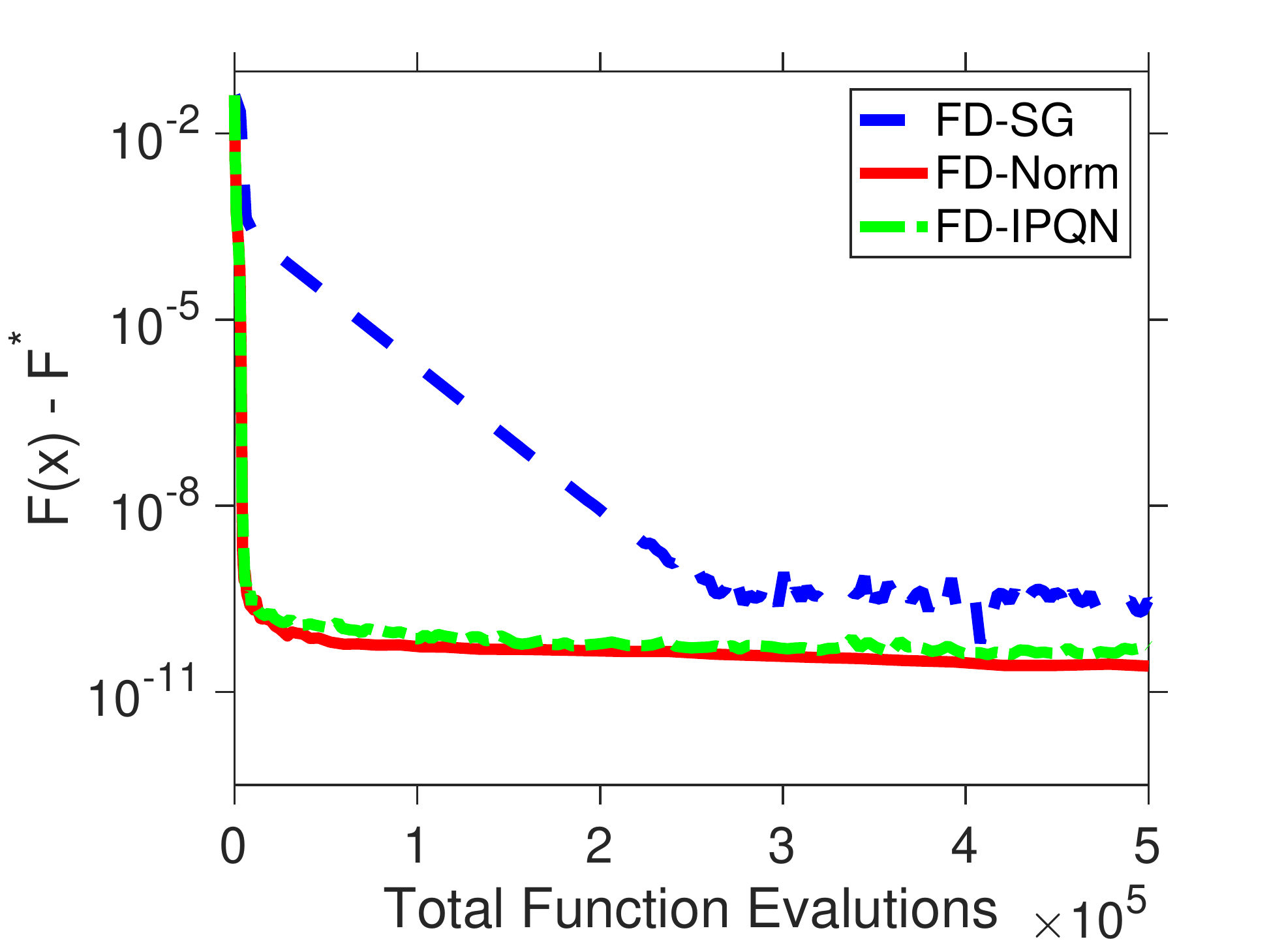} 
		\includegraphics[width=0.40\linewidth]{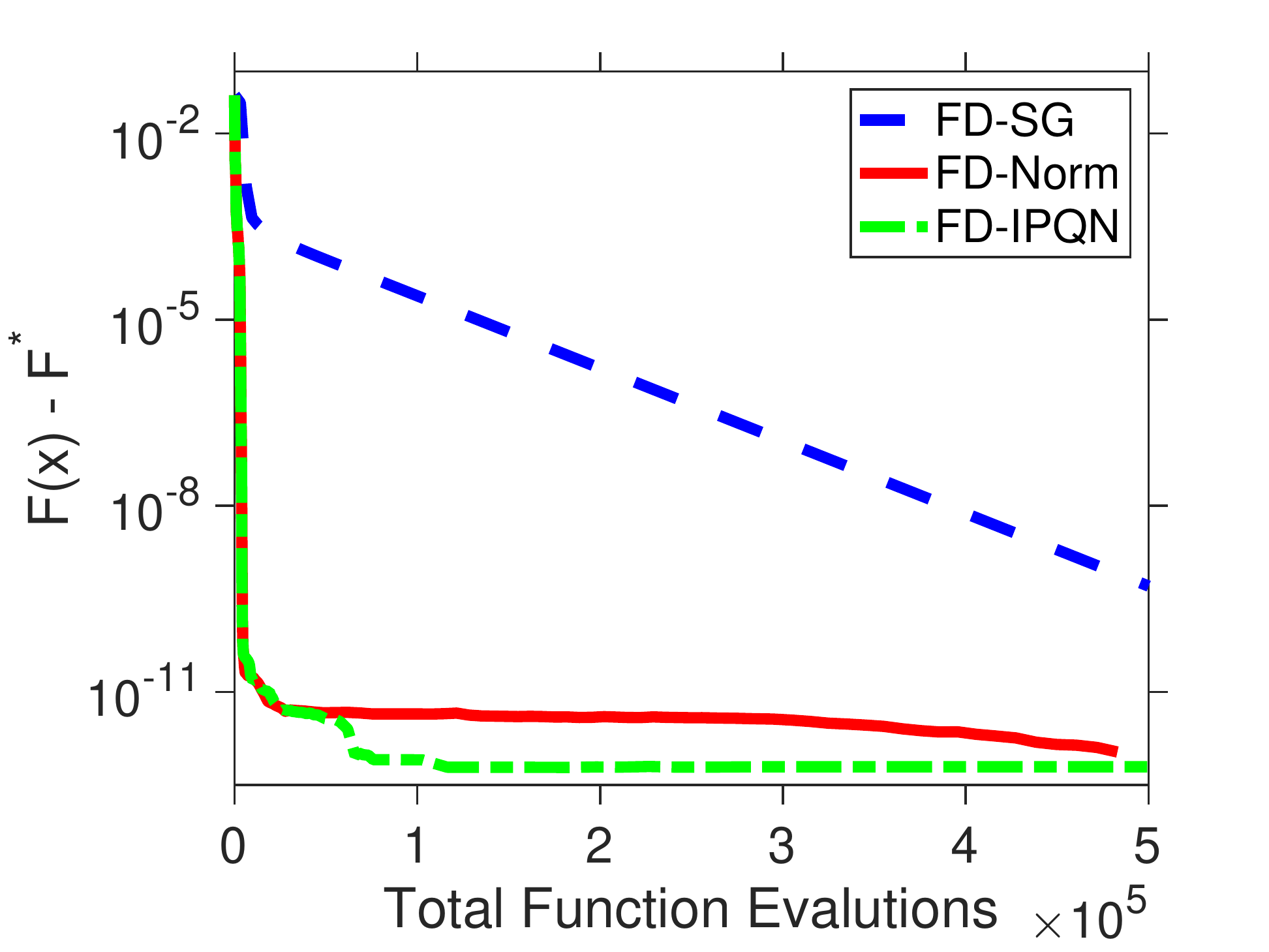}
		\par\end{centering}	
	\caption{Chebyquad function results for a single run: 
		Using $f_{\rm abs}$ (left column) and $f_{\rm rel}$ (right column) with $\sigma=10^{-3}$ (top row) and $\sigma=10^{-5}$ (bottom row).
		\label{fig:Expt1_15}}
\end{figure}

Figure~\ref{fig:Expt1_15} 
\longonly{reports results on the problem $15$ in CUTEr dataset \cite{?} with abs-normal noise and rel-normal noise for two different $\sigma$ values $10^{-3}$ and $10^{-5}$. The vertical axis}measures the error in the function, $F(x) - F^*$,
\longonly{ and the horizontal axis measures }in terms of the total (i.e., including those in the gradient estimates, curvature pair updates, and line search) number of evaluations of $f(x,\zeta)$. The results show that both variants of our finite-difference quasi-Newton method are more efficient than the tuned finite-difference stochastic gradient method. Furthermore, the stochastic gradient method converged to a significantly larger neighborhood of the solution as compared with the quasi-Newton variants.
\longonly{ in the high variance problems $\sigma=10^{-3}$. }
No significant difference in performance was observed between the norm test and inner product quasi-Newton test. These preliminary numerical results show that the modified tests have potential for stochastic problems where the CRN approach is feasible. 
\longonly{Results for the other problems are given in the supplementary material. }

\longonly{\section{Final Remarks}
	We presented finite-difference quasi-Newton methods for solving derivative-free stochastic optimization problems where the sample sizes used in finite-difference gradient estimators are controlled by a modified norm test or inner product quasi-Newton test. These preliminary numerical results show that the modified tests show potential for stochastic problems where the CRN approach is feasible. }

\longonly{In this work, we considered forward finite-differences in all the coordinate directions to estimate the gradients. It is interesting to consider other derivative-free techniques that estimate the gradients in smaller subspaces $(<d)$ which might result in lower computational effort. However, these approaches are challenging and require special attention towards curvature information used in quasi-Newton updates. Another research direction is to consider the discrete reward policy optimization problems which pose significant challenges as the stochastic realizations of the functions are not continuous. We leave these studies as a subject for future research. }
\newpage
\subsubsection*{Acknowledgments}
	This material was based upon work supported by the U.S.\ Department of
	Energy, Office of Science, Office of Advanced Scientific Computing
	Research, applied mathematics and SciDAC programs under Contract No.\
	DE-AC02-06CH11357. 
%




\small

 \bibliography{../../conlearning/bib/ConLearning.bib,../../conlearning/bib/smw-bigrefs.bib,references.bib}

\bibliographystyle{abbrvnat}

\normalsize

\section*{Supplementary Material}

\paragraph{Finite-difference approximation error.}
For any $\nu > 0$ and any function $f(\cdot,\zeta)$ with $L$-Lipschitz continuous gradients, 
we have that
\begin{align}
\|\nabla F^{\FD}_{S_k}(x_k) - \nabla F_{S_k}(x_k)\|^2 &= \sum_{j=1}^{d}\left(\frac{1}{|S_k|}\sum_{\zeta_i \in S_k} \left( \frac{f(x_k + \nu e_j, \zeta_i) - f(x_k,\zeta_i)}{\nu} - \left[\nabla_x f(x_k,\zeta_i)\right]_j\right) \right)^2 \nonumber \\
&\leq \sum_{j=1}^{d}\frac{1}{|S_k|}\sum_{\zeta_i \in S_k}\left(\frac{f(x_k + \nu e_j, \zeta_i) - f(x_k,\zeta_i)}{\nu} - \left[\nabla_x f(x_k,\zeta_i)\right]_j\right)^2 \nonumber \\
&\leq \sum_{j=1}^{d}\frac{1}{|S_k|}\sum_{\zeta_i \in S_k} \left(\frac{L\nu}{2}\right)^2  
= \left(\frac{L\nu}{2}\right)^2d , \nonumber 
\end{align}
where the first equality is by the definitions of $\nabla F_{S_k}^{\FD}(x_k)$ and $\nabla F_{S_k}(x_k) = \frac{1}{|S_k|} \sum_{\zeta_i \in S_k} \nabla_x f(x,\zeta_i)$. The first inequality is due to the fact that $(a_1 + a_2 + \cdots + a_n )^2 \leq n (a_1^2 + a_2^2 + \cdots + a_n^2)$, and the second inequality is due to
\begin{equation*}
f(y, \zeta) \leq f(x,\zeta) + (y-x)^T\nabla_x f(x,\zeta) + \frac{L}{2}\|y-x\|^2, \quad \forall x, y \in \R^d.
\end{equation*}
In a similar manner, we can show that, for the true gradients, the error due to finite-difference approximation is bounded. That is, 
\begin{align}
\|\nabla F^{\FD}(x_k) - \nabla F(x_k)\|^2 &= \sum_{j=1}^{d}\left(\frac{F(x + \nu e_j) - F(x)}{\nu} - \left[\nabla F(x_k)\right]_j\right)^2 \nonumber \\
&\leq \left(\frac{L\nu}{2}\right)^2d \label{eq:ffd_bound}.
\end{align}

\paragraph{Inner product quasi-Newton expected condition.} Consider the following:
\begin{align}
\E_{S_k}[(H_k \nabla F_{S_k}^{\FD}(x_k))^T(H_k \nabla F(x_k))] &= (H_k \nabla F^{\FD}(x_k))^T(H_k \nabla F(x_k)) \nonumber \\
&= (H_k (\nabla F^{\FD}(x_k) - \nabla F(x_k)))^T(H_k \nabla F(x_k)) + \|H_k\nabla F(x_k)\|^2, \nonumber 
\end{align}
where we used
\begin{equation*}
\E_{S_k}[\nabla F^{\FD}_{S_k}(x_k)] =\E_{S_k}\left[\frac{1}{|S_k|}\sum_{\zeta_i \in S_k}\left[\frac{f(x_k + \nu e_j, \zeta_i) - f(x_k, \zeta_i)}{\nu}\right]_{j=1}^{d}\right] 
= \nabla F^{\FD}(x_k),
\end{equation*}
and we have,  
\begin{align*}
(H_k (\nabla F^{\FD}(x_k) &- \nabla F(x_k)))^T(H_k \nabla F(x_k)) \\
&\leq \|H_k\nabla F(x_k)\|\|H_k\|\|\nabla F^{\FD}(x_k) - \nabla F(x_k)\| \\
&\leq \|H_k\nabla F(x_k)\|\|H_k\|\frac{L\nu \sqrt{d}}{2},
\end{align*}
where the last inequality is due to \eqref{eq:ffd_bound}. 
Therefore, 
\begin{align}
\E_{S_k}[(H_k \nabla F_{S_k}^{\FD}(x_k))^T(H_k \nabla F(x_k))] &\geq \|H_k\nabla F(x_k)\|^2 - \|H_k\nabla F(x_k)\|\|H_k\|\frac{L\nu \sqrt{d}}{2} \nonumber \\
&\geq \|H_k\nabla F(x_k)\|\left(\|H_k\nabla F(x_k)\| - \|H_k\|\frac{L\nu \sqrt{d}}{2}\right). \nonumber  
\end{align}
For any $x_k$ such that $\|\nabla F(x_k)\| > \frac{\lambda_{\max}(H_k) L\nu\sqrt{d}}{2\lambda_{\min}(H_k)}$, where $\lambda_{\max}(H_k), \lambda_{\min}(H_k) > 0$ are the largest and smallest eigenvalues of $H_k$,  respectively, it  thus follows that
\begin{align*}
\E_{S_k}\left[(H_k \nabla F_{S_k}^{\FD}(x_k))^T(H_k \nabla F(x_k))\right] > 0. 
\end{align*} 

\paragraph{Numerical experiments setup.}
In the tests of the proposed algorithm, we use $\theta = 0.9$, finite-difference parameter 
$\nu =10^{-8}$, 
L-BFGS memory parameter $m=10$, and line search parameters ($c_1 = 10^{-4}$, $\tau = 0.5$). None of these parameters have been tuned to the problems being considered.


\vspace{3em}

\small
\framebox{\parbox{\linewidth}{
The submitted manuscript has been created by UChicago Argonne, LLC, Operator of 
Argonne National Laboratory (``Argonne''). Argonne, a U.S.\ Department of 
Energy Office of Science laboratory, is operated under Contract No.\ 
DE-AC02-06CH11357. 
The U.S.\ Government retains for itself, and others acting on its behalf, a 
paid-up nonexclusive, irrevocable worldwide license in said article to 
reproduce, prepare derivative works, distribute copies to the public, and 
perform publicly and display publicly, by or on behalf of the Government.  The 
Department of Energy will provide public access to these results of federally 
sponsored research in accordance with the DOE Public Access Plan. 
http://energy.gov/downloads/doe-public-access-plan.}}

\end{document}